\documentclass[12pt]{article}
 
\usepackage{amsmath} 
\usepackage{amsthm} 
\usepackage{amsfonts} 
\usepackage{amssymb}
\usepackage{setspace}
\usepackage{fancyhdr}
\usepackage[dvips]{graphicx}

\newtheorem{thm}{Theorem}[section] 
\newtheorem{defi}[thm]{Definition} 
\newtheorem{prop}[thm]{Proposition} 
 
\newtheorem{cor}[thm]{Corollary}

\newcommand{\R}{\mathcal R} 
\newcommand{\M}{\mathcal M}

\begin{document} 
 
\title{\bf A recipe theorem for the topological Tutte polynomial of Bollob\'as and 
Riordan}

\author{Joanna A. Ellis-Monaghan$^*$ $^1$ 
\and
Irasema Sarmiento $^{**}$}

\maketitle   
\noindent


\footnote{ {\footnotesize $^*$ Department of Mathematics, Saint Michael's College, 1 Winooski Park, Colchester, VT 05439.  
jellis-monaghan@smcvt.edu}  

{\footnotesize $^{**}$ Department of Mathematics, Dipartimento di Matematica, Universit\`a di Roma ``Tor Vergata", Via della Ricerca Scientifica, I-00133, Rome, Italy. sarmient@mat.uniroma2.it}  
  
{\footnotesize $^1$ Support was provided by the National Security Agency and by the Vermont Genetics Network through Grant Number P20 RR16462 from the INBRE Program of the National Center for Research Resources (NCRR), a component of the National Institutes of Health (NIH).  This paper's contents are solely the responsibility of the authors and do not necessarily represent the official views of NCRR or NIH. }  }



\begin{abstract}

 In \cite{BR01}, \cite{BR02}, Bollob\'as and Riordan generalized the 
 classical Tutte polynomial to graphs cellularly embedded in surfaces, i.e. ribbon graphs, thus encoding topological 
 information not captured by the classical Tutte polynomial. We provide a 
 `recipe theorem' for their new topological Tutte polynomial, $R(G)$. We then relate $R(G)$ to the generalized transition polynomial $Q(G)$ of \cite{E-MS02} via a medial graph construction, thus extending the relation between the classical Tutte polynomial and the Martin, or circuit partition, polynomial to ribbon graphs.  We use this relation to prove a duality property for $R(G)$ that holds for both oriented and unoriented ribbon graphs.  We conclude by placing the results of Chumutov and Pak~\cite{CP07} for virtual links in the context of the relation between $R(G)$ and $Q(R)$.

  
\noindent

\end{abstract}


\section{Introduction}\label{introduction}

 
One of Thomas Brylawski's major contributions to the study of the Tutte polynomial was the development of what has come to be known as the `recipe theorem'. It shows that any Tutte-Gr\"{o}thendieck invariant must be an evaluation of the Tutte polynomial, with the necessary substitutions given by the recipe.  This idea first appears in Brylawski's thesis~\cite{Bry70}, with applications throughout much of his early work~\cite{Bry71, Bry72a, Bry72b}.  Overviews and applications of the recipe theorem can be found in his work~\cite{Bry82}, the comprehensive compilation by Brylawski and Oxley~\cite{BO92}, and also in Oxley and Welsh~\cite{OW79}, Welsh~\cite{Wel93}, and the survey by Ellis-Monaghan and Merino~\cite{E-MMa, E-MMb}.

The recipe theorem is essentially a universality statement, and as such, is a very valuable theoretical tool.  Because of this, analogous results are sought for various generalizations of the Tutte polynomial (as well as other graph and matroid polynomials).  Here we find a recipe theorem for a generalization of the Tutte polynomial given by Bollob\'as and Riordan.   

In \cite{BR01}, \cite{BR02}, Bollob\'as and Riordan extended the classical 
Tutte polynomial  to topological graphs, that is, graphs embedded in 
surfaces. In \cite{BR01}, Bollob\'as and Riordan defined the 
{\em cyclic graph polynomial}, a three variable contraction-deletion 
polynomial for graphs embedded in oriented surfaces. They furthered this 
work, using a different approach, in \cite{BR02}, with the 
{\em ribbon graph polynomial}, a four variable 
polynomial for graphs embedded in arbitrary surfaces that subsumes the three variable version. The ribbon graph polynomial, $R(G;x,y,z,w)$, is also sometimes called the
  \emph{Bollob\'{a}s-Riordan polynomial} or 
  \emph{topological Tutte polynomial}.

We provide a `recipe theorem' that, analogously to that for the classical Tutte polynomial, establishes 
conditions for when a graph invariant can be 
calculated from the topological Tutte polynomial and 
gives a formula for this translation. It also restates of the 
universality property of $R(G;x,y,z,w)$ given by Bollob\'as and Riordan~\cite{BR02}.

We show that if certain relations among the 
variables are satisfied, then the topological Tutte polynomial is 
related via an embedded medial graph to the 
generalized transition polynomial of 
\cite{E-MS02}. 
This result extends the relation between the classical Tutte polynomial 
and the Martin polynomial given by Martin~\cite{Mar77} (cf. Jaeger~\cite{Jae90} and Las Vergnas~\cite{Las78, Las79, Las83}). It also 
extends the relation between the Tutte polynomial and the original transition  
polynomial given by Jaeger~\cite{Jae90}. 

We then use these results to give a duality property of $R(G;x,y,z,w)$ and applications to knots and links.

\section{Preliminaries} 
 
We assume the reader is familiar with the work of Bollob\'as and Riordan in ~\cite{BR01}, \cite{BR02} and we adopt the terminology therein, with the conventions of \cite{BR02} taking precedence. We also assume the reader is familiar with cellular embeddings of graphs and with ribbon graphs (also known as fat graphs or band decompositions), and we generally follow Gross and Tucker~\cite{GT87}.  Thus, we only briefly review a few essential concepts.

 A cellular embedding 
  of a graph in a surface (orientable or  
  unorientable) can be specified by providing a sign for each 
  edge and a rotation scheme for the set of
  half edges at each vertex, where  
  a rotation scheme is simply a cyclic ordering of the half edges about a 
  vertex.  This is equivalent to a {\em ribbon graph},
  which is a  
  surface with boundary where the vertices are represented by a set of 
  disks and the edges by ribbons, giving a half-twist to the ribbon of an edge with a 
  negative sign. 
A ribbon graph can also be thought of as a \emph{fat graph}, that is, a 
 slight `fattening' of the edges of the graph as it is embedded in the 
 surface, or equivalently  a `cutting out' of the graph, together with a small
neighborhood of it, from the surface.

\begin{figure}[hbtp] 
     \begin{center} 
               \includegraphics{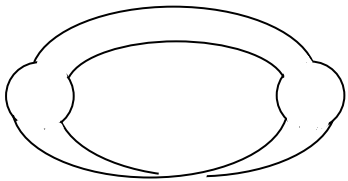} 
       \caption{A ribbon graph which is a M\"obius band with 
         boundary, which can be viewed as a neighborhood of the graph as embedded in a Klein bottle.}\label{Fig:fatgraph}   
     \end{center} 
   \end{figure}

Fig.~\ref{Fig:fatgraph} shows a graph with 
  two vertices and two  
  parallel edges, one positive and one negative.  It is embedded on a 
  Klein bottle, and the ribbon graph is a M\"obius band with boundary.  
 
  For a ribbon graph $G$, we let $k(G), r(G)$, and $n(G)$ be, respectively, the number of connected components, rank, and nullity of the underlying abstract graph.  Additionally, $bc\left( G \right)$ is the number of boundary components of  
  the surface defining the ribbon graph $G$, and $t(G)$ is an index of the
  orientability of the  
  surface, with $t(G)=0$ if the surface is
  orientable, and $t(G)=1$ if it is not.  When $A \in E(G)$, then $r(A)$, $\kappa(A)$, $n(A)$,  $bc(A)$ and $t(A)$, each refer to the spanning subgraph of $G$ on the
  edge set $A$, with embedding inherited from $G$. 

The result of deleting an edge from a ribbon graph is clear.  For contraction of a non-loop edge $e$, assume the sign of $e$ is positive, by flipping one 
  endpoint if necessary to remove the half twist (this reverses 
  the cyclic order of the half edges at that vertex and toggles their
  signs). Then $G/e$ is formed by deleting $e$ and identifying its endpoints into a single vertex $v$. The cyclic order of half edges at $v$ follows first the original cyclic order 
  at one endpoint, beginning where $e$ had been, and continuing with the 
  cyclic order at the other endpoint, again beginning where $e$ had 
  been. The surface that results from sewing disks to the boundaries of a ribbon graph $G$ may not be the same surface that result from sewing disks to the boundary of $G-e$ or $G/e$, particularly when $e$ is a bridge or a loop, i.e. $G-e$ or $G/e$ are not necessarily embedded in the same surface as $G$.

   

There are two definitions of the topological Tutte polynomial, a generating function formulation and a linear recursion formulation, that were shown to be equivalent by Bollob\'as and Riordan~\cite{BR02}.  We begin with the generating function formulation.
   
   \begin{defi}
\label{gen funct TT} Let $G$ be a ribbon
  graph.  Then
  \begin{equation*} 
   \begin{split} 
   R(G;x,y,z,w) &= \sum_{A \subseteq E( G)} 
                 (x - 1)^{r( G ) - r( A )} 
                 y^{n(A)} z^{\kappa(A) - bc(A) + n(A)} 
                   w^{t(A)}\\ 
                &\in Z[x,y,z,w]/\vphantom {Z[x,y,z,w]} 
                  \langle w^2  - w\rangle . 
   \end{split} 
  \end{equation*} 
 
  \end{defi}

   The linear recursion formulation derives from the following theorem.

   

\begin{thm}~\cite{BR02}\label{ribbon delete-contract} 
  If $G$ is a ribbon graph, then  
 \begin{equation*} 
R(G;x,y,z,w) = R(G/e;x,y,z,w) + R(G - e;x,y,z,w),
\end{equation*}
if $e$ is an ordinary edge; and
\begin{equation*}
R(G;x,y,z,w) = x\,R(G/e;x,y,z,w),
\end{equation*}  
if $e$ is a bridge. 

\end{thm}

   Repeated application of this theorem reduces a ribbon graph to a 
  disjoint union of embedded \emph{bouquets},
   that is, embedded graphs each  
  consisting of a single vertex with some number of loops, and the topological 
  information is distilled into these minors of the original graph. 

\begin{figure}[hbtp] 
     \begin{center} 
               \includegraphics{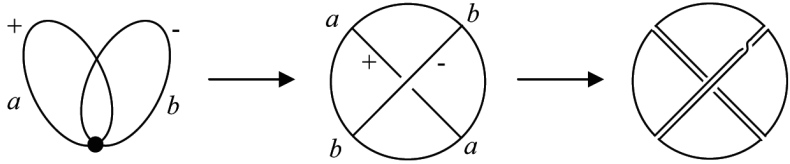} 
       \caption{An embedded bouquet with two loops, one positive and one negative, with its signed chord diagram and the boundary components of the fattened signed chord diagram.}\label{Fig:figureB}   
     \end{center} 
   \end{figure}

Signed chord diagrams are a useful device for determining the  
  relevant parameters of an embedded bouquet.  A
  \emph{signed chord  diagram} is a circle with $n$ symbols  each appearing twice on its
  perimeter, with a signed $(\pm 1)$ chord drawn between each pair of 
  like symbols.  
 We assign a symbol to each loop of an embedded bouquet 
  $G$  and arrange them on the 
  perimeter of the circle in the chord diagram in the same order as
  the cyclic order of the half-edges  
  about the vertex, with a chord receiving the same sign as the loop it 
  represents.  Since signed chord diagrams are exactly equivalent to bouquets, we will use the terms interchangeably.  If we `fatten' the chords as in 
  Fig.~\ref{Fig:figureB}, with a negative  
  chord receiving a half-twist, then $bc(G)=bc(D)$, the number of 
  components in the resulting diagram.   Similarly, since $G$ has only one vertex, $n(G)$ is the number of 
  edges of $G$, which is the number of chords of $D$, so we denote this 
  by $n(D)$.  We also set
  $t(D)=t(G)$, and note that
  $t (D) = 0$   if all  
  chords of $D$ have a positive sign, and $t(D)=1$
  otherwise.  This,  
  combined with  Definition~\ref{gen funct TT}, gives in Proposition ~\ref{blossom eval} the 
 necessary evaluations of the terminal forms to complete the linear recursion formulation.

  \begin{prop}\label{blossom eval} 
   If $G$ is an embedded bouquet with 
  corresponding signed chord diagram $D$, then  
  \[ 
   R(G;x,y,z,w) = \sum_{D' \subseteq D} y^{n(D')}   
                                      z^{1- bc(D')+ n(D')} 
                                      w^{t(D')}, 
  \] 
  where the sum is over all subdiagrams $D'$ of $D$.  
  \end{prop}

\section{The recipe theorem}\label{recipethm} 
  
In this section we give the recipe theorem that specifies precisely when and how a function may be recovered from $R(G;x,y,z,w)$. It is essentially a restatement of the universality of $R(G;x,y, z,w)$ from Bollob\'as and Riordan~\cite{BR02} in a form that facilitates its application.  
  
\begin{figure}[hbtp] 
     \begin{flushleft} 
               \includegraphics{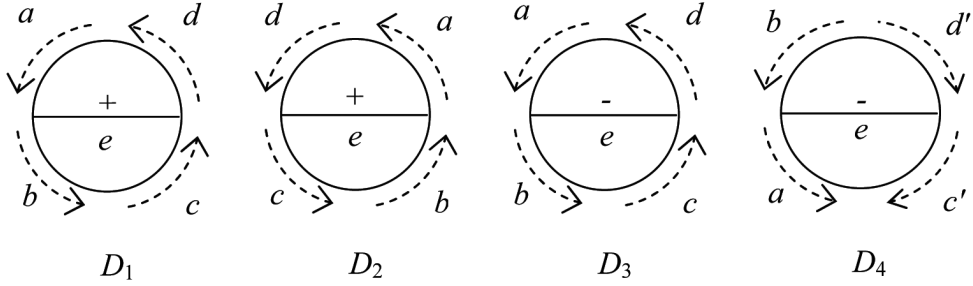} 
       \caption{Related chord diagrams where the dotted arrows indicate directions of sequences of labels along the perimeter.}   \label{Fig:relateddiagrams}   
     \end{flushleft} 
\end{figure}

Following \cite{BR02} we say that two chord diagrams are related by a
\emph{rotation} about the chord $e$ if they are related as
$D_1$ and $D_2$ in Figure~\ref{Fig:relateddiagrams}. Two chord diagrams are related by a
\emph {twist} about  $e$ if they are related as
$D_3$, $D_4$ in Fig.~\ref{Fig:relateddiagrams}, where a letter represents a sequence of labels about the circle, and a prime symbol means to reverse the order of the sequence.  Two diagrams are \emph{related} if they
are related by a sequence of rotations and twists. From \cite{BR02} we have that any diagram is
related to a \emph{canonical signed chord diagram} $D_{ijk}$ $(0\leq k\leq 2$) that consists
of $i-2j-k$ positive chords intersecting no other chord, $j$ pairs
of intersecting positive chords and $k$ negative chords intersecting
no other chord.

\noindent  

We also observe that 
 \begin{align} \label{loops}
R(D_{100};x,y,z,w)=1+y, \\ 
R(D_{101};x,y,z,w)=1+yzw,\\ 
\textnormal{ and } R(D_{210};x,y,z,w)=y^2z^2+2y+1.
\end{align} 

Since from \cite{BR02},  $R(G;x,y,z,w)$ is multiplicative on one-point joins of ribbon graphs, where the new rotation is given by simply concatenating the rotation systems of the vertices being joined, then 
\begin{equation} \label{Dijk}
R(D_{ijk})=[R(D_{100})]^{i-2j-k}[R(D_{210})]^j[R(D_{101})]^k.  
\end{equation}

Furthermore, Bollob\'as and Riordan~\cite{BR02} show that if $D$ is equivalent to $D_{ijk}$ then $R(D)=R(D_{ijk})$.

\begin{thm}\label{recipe}(The recipe theorem.) 

Let $F$ be a map from a minor closed subset ${\mathcal M}$ of ribbon 
graphs containing $D_{210}, D_{100}$, and $D_{101}$ to a commutative ring $\R$ with 
unity; let $s=F (D_{210})$, $q=F (D_{100})$ and  
$ r=F (D_{101})$;  and suppose there are elements $\alpha, x,u,v \in R$ with $\alpha$ a unit such that: 
 
\begin{enumerate}  
  
\item  \label{delcon}
  
$$  
F (G)=  
\begin{cases}  
F(G/e)+F(G\setminus e){\text{ if }} e   
{\text{ is  ordinary }},\\  
xF(G/e) {\text{ if }} e {\text{ is a bridge}}; 
\end{cases}  
$$  
  
\item \label{blockmult}
 
$F(G\dot{\cup}H )= F (G)F(H)$ and $\alpha F(G*H)=  F (G)F(H)$ where $G$, $H$ are embedded bouquets, $G\dot{\cup}H$ is the disjoint union, and $G*H$ is the one-point join, again with concatenated rotation system;

\item 
 
$F (\bold{E})=\alpha^n$ if $\bold{E}$ is an edgeless graph on $n$ vertices;

\item \label{relations}

$(q-\alpha)^2u^2=\alpha (s-2q+\alpha)$, and $(q-\alpha)uv =r-\alpha $, and also $v=v^2$. 
  
\end{enumerate}  
 
Then 
 
$$F (G)=\alpha^{k(G)}R( G ;x,\alpha ^{-1}q-1, u,v),$$
 
where $k(G)$ is the number of components of $G$. 
 
\end{thm}

\begin{proof} 
The proof is by a double induction, first on the number of chords in a signed chord diagram, and then on the number of non-loop edges of $G$.

We first note that by Item~\ref{blockmult} and~\ref{relations} and Equations ~\ref{loops} through~\ref{Dijk}, the result holds for any canonical signed chord diagram.

Since this recipe theorem is also a universality statement, unsurprisingly the proof uses the same central observations about chord diagrams as the proof of universality for Theorem $2$ from Bollob\'as and Riordan~\cite{BR02}. Thus,  from \cite{BR02}, we have that since $F$ satisfies Item \ref{delcon}, it satisfies
\begin{eqnarray}\label{Frelated}
F(D_1)-\mu F(D_1^{'})=F(D_2)-\mu F(D_2^{'})\label{Frelated1}
\textnormal{ and}
\\
F(D_3)-\mu F(D_3^{'})=F(D_4)-\mu F(D_4^{'}),\label{Frelated2}
\end{eqnarray}

\noindent where the $D_i$'s are related as in Figure~\ref{Fig:relateddiagrams} with $D_i^{'}=D_i-e$, and $\mu =1$ if there is a chord from $a \cup b$ to $c \cup d$, and otherwise $\mu =x$.

The same identities hold for $R(G;x,y,z,w)$, and hence hold for $F'=F-\alpha R$. 
 
Now suppose $G=D$ is a signed chord diagram corresponding to an embedded bouquet. 
Assume by induction that 
$F (D)=\alpha R(D;x,\alpha ^{-1}q-1,u,v)$ if $D$ is a signed chord diagram with 
fewer than $m$ chords. Then $F'(G)$ vanishes on signed 
chord diagrams with fewer than $m$ chords.

This, with Equations~\ref{Frelated1} and ~\ref{Frelated2}, implies that
$F'(D)=F'(D')$ if $D$ and $D'$ are related
chord diagrams with $m$ chords. In particular,  $D$ is related
to a
canonical diagram $D_{ijk}$, so $F'(D)=F'(D_{ijk})$. 
 
Thus 
 
$$F'(D)=F'(D_{ijk})=F(D_{ijk})-\alpha R(D_{ijk};x,\alpha ^{-1}q-1,u,v)=0,$$

\noindent
and hence by induction, $F'(D)=0$ on all signed chord diagrams. This extends to disjoint unions of embedded bouquets by Item~\ref{blockmult} and that $R(G;x,y,z,w)$ is also multiplicative on disjoint unions.
 
Thus the result holds for all ribbon graphs with no non-loop edges. If $G$ has a non-loop edge $e$, the result is immediate by induction from Item~\ref{delcon}, observing that, if $e$ is a bridge, then $k(G/e)=k(G)$.

\end{proof} 


In analogy to the classical case, we call a function on ribbon graphs satisfying the conditions of Theorem~\ref{recipe} a \emph{topological Tutte invariant}, and the theorem itself justifies calling the ribbon graph polynomial of Bollob\'as and Riordan \emph{the} topological Tutte polynomial.

We give a quick example by applying Theorem~\ref{recipe} to give   
the relationship noted 
in \cite{BR02} between $R(G;x,y,z,w)$ and 
the oriented graph invariant $C(G;x,y,z)$ of \cite{BR01}. Cyclic graphs form a minor-closed subset of ribbon graphs, and we 
extend $C$ very slightly by defining $C(D_{101};x,y,z)=1+yz^{\frac{1}{2}}w$, noting that the domain is still 
minor-closed. If we let $\R$ be the quotient ring
${\mathbb Z}[x,y,z^{\frac{1}{2}},w]/\langle w-w^2 \rangle $, then $C$ satisfies the recipe theorem with 
$\alpha=1$, $1+2y+y^2 z=s$, $1+y=q$, 
$1+yz^{\frac{1}{2}}w=r$ and taking $u=z^{\frac{1}{2}}$ and
$v=w$. Thus, $C(G;x,y,z)=R(G;x,y,z^{\frac{1}{2}},w)$. 


The polynomial $R$ has the property that it immediately identifies whether or not a ribbon graph is oriented, just by the absence or presence of the one variable $w$.  For an arbitrary topological Tutte invariant however, this property is sensitive to the structure of the ring in which it takes its values. 
 
 
\begin{cor}\label{discern} 
 
If $F$, $\M$, $\R$ satisfy Theorem~\ref{recipe}, with both $q-\alpha$ and $r-\alpha$ being units of $\R$, then $v$=1, and thus $F$ does not discern orientation by the presence or absence of a single idempotent element. 
 
\end{cor} 
 
\begin{proof} 
 
From Item~\ref{relations}, we have that  
$\frac{r-\alpha}{q-\alpha}=uv=uv^2=v \frac{r-\alpha}{q-\alpha}$. Thus, since $q-\alpha$ and $r-\alpha$ are units,  $v=1$,  and hence 
$F (G)=\alpha^{k(G)}R(G;x,\alpha ^{-1} q-1,u,1)$.
 
\end{proof} 

This corollary raises a number of questions.  Suppose $F$ is a topological Tutte invariant with the properties of Corollary~\ref{discern}. Then $F$ does not determine orientation by the presence or absence of one idempotent element.  However, unless $u$ also equals 1, $F$ can still distinguish between oriented and unoriented embeddings of a graph. For example, it is easy to check that $R(G,x,y,z,1)$ distinguishes all the canonical chord diagrams, and furthermore a chord diagram $D$ is orientable if and only if a term of the form $yz$ does not appear in $R(D;x,y,z,1)$.  Since this is so, is $w$ strictly necessary to record orientability information?  I.e. $F$ can obviously be computed from $R$, but is it possible to recover $R$ from some such $F$?   We suspect not, since a consistent translation from $F$ to $R$ is problematic even on the canonical chord diagrams, which leads to the next question.  Is $R$ actually  more refined than any such $F$?  I.e. is there a pair of ribbon graphs distinguished by $R$ that are not distinguished by $F$?  Finally, the most basic question is whether it is always possible to determine from some such $F$ that a ribbon graph $G$ is oriented.

\section{Medial graphs and transition polynomials}\label{medials} 
 
The classical Tutte polynomial, $T(G;x,y)$, among many other properties, encodes information about families of Eulerian circuits in the medial graph of a planar graph.  This theory is the result of a relation between the classical Tutte polynomial and the Martin, or circuit partition, polynomial.  Here we extend this theory to ribbon graphs, giving an analogous result relating the topological Tutte polynomial of a ribbon graph to the transition polynomial of its topological medial graph, where the transition polynomial of \cite{E-MS02} is a multivariable generalization of the circuit partition polynomial.   The original relation between the Tutte polynomial and the Martin polynomial can be found in Martin's 1977 thesis \cite{Mar77}, with the theory considerably extended by Martin~\cite{Mar78}, Las Vernas~\cite{Las79, Las83, Las88}, Jaeger~\cite{Jae90}, Bollob\'as~\cite{Bol02}, and \cite {E-M00, E-M04a, E-M04b, E-MS02}.  An overview can be found in Ellis-Monaghan and Merino~\cite{E-MMa, E-MMb}.

The \emph{medial graph} of a
 connected  planar graph $G$ 
 is constructed by placing a vertex on each edge of $G$ and drawing
 edges around the faces of $G$.  The faces of this medial graph are
 colored black or white, depending on whether they contain or do not
 contain, respectively, a vertex of the original graph $G$.  This face
 two-colors the medial graph.  The edges of the medial graph are then
 directed so that a black face is to the left of each edge.  This directed medial graph is denoted $\vec{G}_m$, and is an \emph{Eulerian digraph}, that is, the number of incoming edges is equal to the number of  outgoing edges at each vertex.

For the circuit partition polynomial we first recall that an \emph{Eulerian vertex state}  is a choice of reconfiguration at a vertex of
an Eulerian digraph $\vec{G}$. The reconfiguration consists of replacing a $2n$-valent vertex $v$
with $n$ 2-valent vertices joining pairs of edges originally adjacent to $v$, where each incoming edge must be paired with an outgoing edge.  An \emph{Eulerian graph state} of an Eulerian digraph $\vec{G}$ is the result of choosing
one vertex state at each vertex of $\vec{G}$. Note that a graph state is a disjoint union of consistently oriented
cycles. Let $\cal{S}$ denote the set of Eulerian graph states of an Eulerian digraph $\vec{G}$, where $\cal{S}$ is not "up to isomorphism", so that each individual state is included in the set.  

The \emph{circuit partition polynomial} of an Eulerian digraph is $j\left( {\vec G}; x \right) = \sum\limits_{S \in \cal{S}} {x^{c(S)} }$, where $c(S)$ is the number of components of the state $S$.  For a connected planar graph $G$ with oriented medial graph $\vec {G}_m$, the relations among the Martin polynomial, circuit partition polynomial, and classical Tutte polynomial are

\begin{equation}\label{martutte}
x^{k(G)} m(\vec {G}_m ;x+1) = j(\vec {G};x)=  x^{k(G)}T(G;x+1,x+1).
\end{equation}

  In the context of transition polynomials such as the circuit partition polynomial (and also certain link invariants), the number of components of a state does not count isolated vertices.  Thus, we use $c$ here for components, in contrast with the $k$ we use in the context of Tutte polynomials where isolated vertices are included in the component count.

\begin{figure}[hbtp] 
     \begin{center} 
               \includegraphics{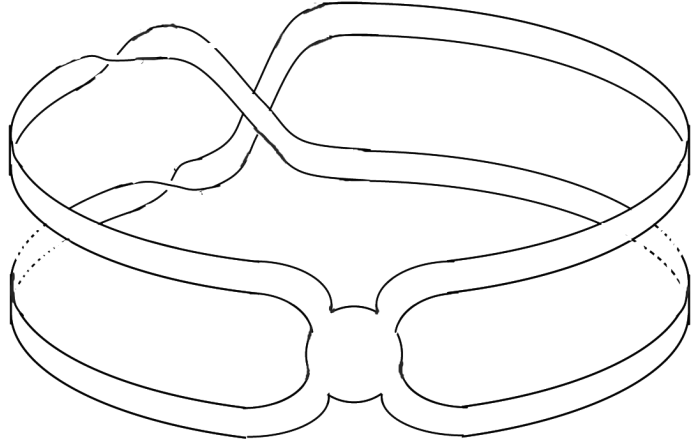} 
       \caption{The topological medial graph of a single loop with a negative edge.}
   \label{Fig:medial}   
     \end{center} 
\end{figure}

To extend Equation (\ref{martutte}) to ribbon graphs, we begin with the notion of a topological medial graph.   Let $G$ be a connected ribbon graph, thought of as being cellularly embedded in a surface.  We construct the medial graph $G_m$ in the same surface, exactly as in the planar case.  That is, we place a vertex of $G_m$ on each edge of $G$, and draw the edges of $G_m$ by following two adjacent half-edges of $G$ around the face they bound, as in Figure~\ref{Fig:medial}, which shows the topological medial graph $G_m$, where $G$ is a single loop with a negative edge. Both $G_m$ and $G$ are embedded in a Klein bottle.  Note that if $e$ is a negative edge, two of the half edges which are consecutive with the $e$ between them receive a half twist.  Whether or not an edge of $G_m$ is positive or negative depends on the parity of the half twists on its half edges.

\noindent 

We define the medial graph of an isolated vertex to be a ``free loop'', 
that is an edge, but no vertex, following the boundary of a small disk on 
the surface containing the vertex.

The generalized transition polynomial, $Q(G; W,t)$  of \cite{E-MS02} is a multivariable extension of the circuit partition polynomial that assimilates the   transition polynomial of Jaeger~\cite{Jae90} for 4-regular graphs and generalizes it to arbitrary Eulerian graphs. We recall the essentials needed for the current application, and refer the reader to \cite{E-MS02} for full details.

A \emph{weight system}, $W( G )$, of
    an Eulerian graph $G$ is an assignment of a pair weight in a unitary ring $\R$ to every
    possible pair of adjacent half edges of $G$. (We simply write $W$ for $W(G)$ when the graph is clear from context.) A \emph{pair weight} is the particular
    value $p\left( {e_v ,e'_v } \right)$ associated by the weight system
    to a pair of half edges $e_v ,e'_v$ incident with a vertex $v$. The \emph{vertex state weight}   of a
    vertex state is $\prod {p( {e_v ,e'_v } )}$
    where the product is over the pairs of half edges comprising
    the vertex state. The \emph{state weight} of a graph state
    $S$ of a graph $G$ with weight system $W$ is $\omega( S ) =
    \prod{\omega( {v,S} )}$,
    where $\omega( {v,S} )$ is the vertex state weight of
    the vertex state at $v$ in the graph state $S$, and where the
    product is over all vertices of  $G$.

The generalized transition polynomial is defined exactly like the circuit partition polynomial, simply with the addition of keeping track of the specific weights for each pair of adjacent edges given by the weight system.

\begin{defi}\label{transpolydef} Let $G$ be a graph having weight system $W$ with values in a unitary ring $\R$.  Then the generalized transition polynomial is $Q(G; W,t)= \sum\limits_{S \in \cal{S}} {\omega( S )t^{c(S)} }$, where $c(S)$ is the number of components of the state $S$.  

\end{defi}

In the case that $G$ is a planar graph with oriented medial graph $\vec{G_m}$, we can assign a weight system $W$ to the underlying medial graph $G_m$, with pair weights of 1 for each pair of half edges where one is incoming and the other outgoing in $\vec{G_m}$, and 0 otherwise.  With this weight system, $Q(G_m; W, x)= j(\vec{G_m};x)$.

\begin{figure}[hbtp] 
     \begin{center} 
               \includegraphics{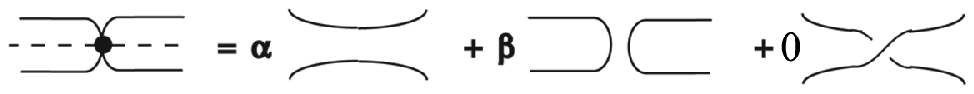} 
       \caption{The weight system for a topological medial graph.  On the left is a neighborhood of a vertex of the topological medial graph $G_m$ with the corresponding edge of $G$ as a dotted line.  On the right are the three possible vertex states, uncut, cut, and crossing, with their respective vertex state weights.}
   \label{Fig:weightsystem}   
     \end{center} 
\end{figure}

For the current application, we will restrict $Q$ to medial ribbon graphs.  In the rotation system about a vertex $v$ of a medial graph $G_m$ of a graph $G$, we can consider \emph{six} half edges, the four half edges actually belonging to $G_m$, plus the two half edges of the edge of $G$ corresponding to $v$.  This allows us to define the following weight system, which we will refer to hereafter as the \emph{medial weight system}.  If a pair of edges of $G_m$ are consecutive in the rotation system at $v$ without a half edge of $G$ between them, we assign them a pair weight of $\sqrt{\alpha}$.  If they are consecutive with a half edge of $G$ between them, we assign them a pair weight of $\sqrt{\beta}$. Otherwise, their pair weight is zero.  The square root is just a notational convenience so that the vertex state weights of the two possible nonzero vertex states will be either $\alpha$ or $\beta$.  This refinement to the level of pair weights is not strictly necessary to the current paper, but we provide it since it is required for splitting formulas of $Q$ in \cite{E-MS02} that may be useful in future applications. We can think of these two vertex states as either leaving a ribbon of $G$ intact, or `snipping' through it, and so we will refer to them as `uncut' or `cut' vertex states, respectively. See Figure~\ref{Fig:weightsystem}.

 \noindent  

We now give the relationship between the generalized transition polynomial and the topological Tutte polynomial of Bollob\'as and Riordan that extends the classical case.  Although the proof technique here is very similar to that of Moffatt~\cite{Mof08} and Chmutov and Pak~\cite{CP07}, the result is much broader, since the link invariants they address are specializations of the generalized transition polynomial.  See Section~\ref{links}.

\begin{thm}\label{transpoly} 

Let $G$ be a ribbon graph with topological medial graph $G_m$, and let $G_m$ have the medial weight system $W$.  Then

\[
Q(G_m; W,t)=\alpha^{r(G)} \beta^{n(G)} t^{k(G)} R \left ( G; \frac {\beta t}{\alpha} +1, \frac {\alpha t}{\beta}, \frac {1}{t},1 \right ). 
\]

\end{thm}

\begin{proof}
Observe that if $W$ is the medial weight system, then $Q(G_m; W, t) = \sum\limits_{S \in \cal{S}} {\alpha^{a(S)} \beta^{b(S)} t^{c(S)} }$, where $a(S)$ and $b(S)$ are the number uncut or cut vertex states, respectively, in the graph state $S$.  Thus, we can use the edges of $G$ to index this sum, and thinking of $H$ as the set of edges that are uncut, we have that 
\[
Q(G_m; W, t) = \sum\limits_{H \subset E(G)} {\beta^{|E|-|H|} \alpha^{|H|} t^{bc(H)} } =\beta^{|E|} \sum\limits_{H \subset E(G)} { \left ( \frac{ \alpha}{\beta} \right ) ^{|H|} t^{bc(H)} }.
\]   

We then note that

\begin{multline*}
R(G;,x,y,z,1)= \\
(x-1)^{-k(G)}y^{-v(G)}z^{-v(G)}\sum\limits_{H \subset E(G)}{ \left ( (x-1)yz^2) \right ) ^{k(H)} (yz) ^{|H|}z^{-bc(H)}}.
\end{multline*}

Substituting $z=\frac{1}{t}$, $y=\frac{\alpha t}{\beta}$ and $x=\frac{\beta t}{\alpha} +1$  yields the result.

\end{proof}

 Alternatively, Theorem~\ref{transpoly} could also have been proved as easily using Theorem~\ref{recipe}.

\section{Duality}\label{duality} 

The results of Section~\ref{medials} provide tools to determine properties of $R(G;x,y,z,w)$.  

We write $G^*$ for the dual ribbon graph of a ribbon graph $G$ (see Gross and Tucker~\cite{GT87}  or Bollob\'as and Riordan~\cite{BR02}), and note that $G_m^*$ and $G_m$ are isomorphic ribbon graphs. 
Furthermore, $G_m$ is a four regular ribbon graph that determines 
the same surface as $G$ and  $G^*$ do.  The cellular embedding requirement implies that all of  $G$, $G^*$, $G_m$, and$G_m^*$ have the same number of components.

\begin{figure}[hbtp] 
     \begin{center} 
               \includegraphics{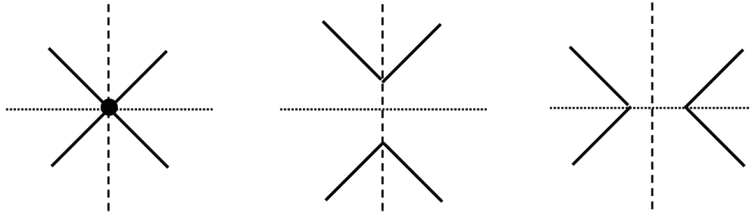} 
       \caption{The lefthand figure shows a vertex and four half edges of $G_m$, or equivalently $G^*_m$, (solid) with the edge of $G$ (dotted) and $G^*$ (dashed).  The center figure shows the vertex state which is uncut with respect to $G$ and cut with respect to $G^*$.  The righthand figure shows the vertex state which is cut with respect to $G$ and uncut with respect to $G^*$.}
   \label{fig:duality}   
     \end{center} 
\end{figure}

If $G_m$ is a topological medial graph with weight system $W$, then the \emph{dual weight system} $W^*$ results from exchanging the roles of $\alpha$ and $\beta$.  This leads to the following duality relation for the generalized transition polynomial, the central idea for which is illustrated in Figure~\ref{fig:duality}.

 \noindent  

\begin{thm}\label{transitiondual}
If $G$ is a ribbon graph with dual $G^*$, then

\[
Q(G_m^*; W(G_m^*), t) = Q(G_m; W^*(G_m),t).
\]
\end{thm}

\begin{proof}

If we think of $H \subset E(G)$ as indexing the uncut edges, then 
\[
Q(G_m; W^*(G_m), t) = \sum\limits_{H \subset E(G)} {\alpha^{|E|-|H|} \beta^{|H|} t^{bc(H)} }.
\]
  However, as in Figure~\ref{fig:duality}, an uncut vertex state of $G_m$ corresponds to a cut vertex state of $G_m^*$ and vice versa.  Also $bc(G|_H) = bc(G^*|_{E-H})$.  Thus, $Q(G_m; W^*(G_m), t) = \sum\limits_{H \subset E(G^*)} {\alpha^{|E|-|H|} \beta^{|H|} t^{bc(E-H)} }$, where we think of $H \subset E(G^*)$ now as indexing the cut edges of $G*$.  But if we instead index over sets of \emph{uncut} edges, this is then $\sum\limits_{H \subset E(G^*)} {\beta^{|E|-|H|} \alpha^{|H|} t^{bc(H)} } = Q(G_m^*; W(G^*_m), t).$
\end{proof}

We can now give a duality relation that extends the duality relation 
for $R(G;x,y,z,w)$ given by Bollob\'as and Riordan in \cite{BR02} from one degree of 
freedom to two, thus giving a natural extension of the duality of the 
classical Tutte polynomial.  This theorem was first announced in \cite{E-MS05} and has since been referenced by Moffatt~\cite{Mof, Mof08}, Chmutov~\cite{Chm}, and Vignes-Tourneret~\cite{V-T}. It is stronger than the version in \cite{Mof08} in that it applies to unoriented as well as oriented ribbon graphs. Chmutov~\cite{Chm} gives an alternative proof and slightly different formulation.

\begin{thm} 
 
Let $G$ be cellularly embedded in a not necessarily connected surface $\Sigma$, let $G^*$ be its dual, and let $\gamma=2 \sum{g_i} +\sum{g_i}$, where the first sum is of the genera of the orientable components of $\Sigma$ and the second sum is of the genera of the unorientable components.  Then 
 
\begin{equation}\label{duals}
\beta ^ {\gamma}R(G^*;\frac{\beta t}{\alpha} +1, \frac{\alpha t}{\beta}, 
 \frac{1}{t}, 1)=
 \alpha ^ {\gamma}R(G;\frac{\alpha t}{\beta} +1, \frac{\beta t}{\alpha}, 
 \frac{1}{t}, 1).
\end{equation}

  Furthermore, if we write  
  $\sqrt{t^2}$ as $t$, and $\sqrt{\frac{\beta^2}{\alpha^2}}$ as $\frac{\beta}{\alpha}$,  then we may substitute $x=\frac{\beta t}{\alpha}$ 
  and $y=\frac{\alpha t}{\beta}$ to rewrite this as 
 
$$ 
x^{\frac{\gamma}{2}} R(G;1+x, y, \frac{1}{\sqrt{xy}}, 1)= 
  y^{\frac{\gamma}{2}}R(G^*;1+y, x, \frac{1}{\sqrt{xy}},1).$$

If $\alpha=\beta=1$, then Equation~(\ref{duals}) reduces to the one variable duality identity given in 
\cite{BR02}. 
\end{thm}

\begin{proof} 
Theorems~\ref{transpoly} and \ref{transitiondual} give that 

\begin{equation*}
\alpha^{r(G^*)} \beta^{n(G^*)} t^{k(G^*)} R(G^*;\frac{\beta t}{\alpha} +1, \frac{\alpha t}{\beta}, 
 \frac{1}{t}, 1)=
 \beta^{r(G)} \alpha^{n(G)} t^{k(G)}R(G;\frac{\alpha t}{\beta} +1, \frac{\beta t}{\alpha}, 
 \frac{1}{t}, 1).
\end{equation*} 

To simplify the exponenets, we use the invariance of the Euler characteristic on each connected component of $\Sigma$, namely that $v_i-e_i+f_i=2-\gamma_i$, where $\gamma_i$ is $2g_i$ or $g_i$ depending orientability, and where $v_i$, $e_i$, and $f_i$ are the number of vertices, edges, and faces, respectively, of $G$ on the $i^{th}$ component.   By definition and duality, $r(G^*)=v(G^*)-k(G^*)=f(G)-k(G)$. Then $f(G)-k(G) =\sum_{i}{f_i}+\sum_{i}{1}$, where the sum is over the number of components of $\Sigma$. By the invariance of the Euler characteristic on each component, this becomes $\sum_{i}{(2-\gamma_i+e_i+v_i-1)}=k(G)+e(G)-v(G) - \gamma$.  Thus, $r(G^*)+\gamma =n(G)$, and similarly,  
$r(G)+\gamma =n(G^*)$, which gives the result.
\end{proof}

\section{Applications to links}\label{links} 

We now turn to the virtual links of Kauffman~\cite{Kau99} and Goussarov, Polyak, and Viro~\cite{GPV00}.  Chmutov and Pak~\cite{CP07} found a relation between the Kauffman bracket of virtual links and $R(G)$. We will focus on their result for signed graphs, since it subsumes their unsigned version.  Here we show that since the Kauffman bracket is another specialization of the generalized transition polynomial, $Q(G)$, the results of \cite{CP07} follow immediately from those of Section~\ref{medials}.

\begin{figure}[hbtp] 
     \begin{center} 
               \includegraphics{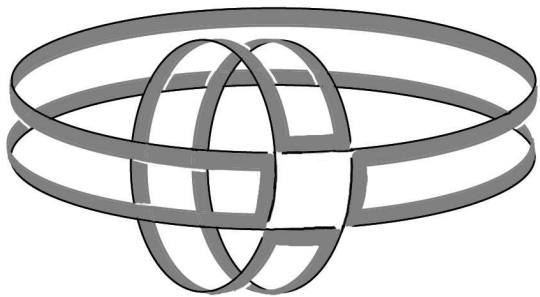} 
       \caption{A checkerboard coloring of a plane drawing of a four component link in the torus. The four crossings forming a square in the center are \emph{classical}, while the rest are \emph{virtual}. }
   \label{fig:checkerlink}   
     \end{center} 
\end{figure}

Let $\Sigma$ be a compact oriented surface.  Here we view a virtual link $L$ as an unoriented link
in $\Sigma\times [0,1]$,
with link diagram $\tilde{L}$ in $\Sigma$, such that the link universe $\Gamma_L$ is cellularly embedded in $\Sigma$ (and hence is a ribbon graph).  In a plane drawing of $L$, the crossings corresponding to actual crossings are called \emph{classical}, and the others, artifacts of the projection, are called \emph{virtual}.  Following Kamada~\cite{Kam02, Kam04}, a plane drawing of a virtual link is \emph{checkerboard colorable} if a small neighborhood of one side of each strand can be colored so that a checkerboard pattern is formed at classical crossings, while the strand coloring passes through virtual crossings unchanged. See Figure~\ref{fig:checkerlink}.  We also recall that an $A$ splitting at a classical crossing is the result of opening a channel between the two regions swept out by rotating the top strand counterclockwise, and a $B$ splitting joins the other two regions.

\noindent 

\begin{defi}\label{bracket} 
The {\em generalized Kauffman bracket} of a virtual link diagram
$\tilde{L}$ is the polynomial 

$$[\tilde{L}](A,B,d)=\sum_{S\in S(\tilde{L})}A^{a(S)}
B^{b(S)}d^{c(S)-1},$$

\noindent
where $S(\tilde{L})$ is the set of states of $\tilde{L}$, where
$a(S)$ and $b(S)$ are the number of
$A$ and $B$ splittings, respectively, in the state $S$ and where
$c(S)$ is the number of components.
\end{defi}

If $\Gamma_L$ is the universe of a virtual link digram $\tilde{L}$, then it inherits a weight system, $W_L$, from $\tilde{L}$ by assigning a pair weight of $\sqrt{A}$  to half edges which are joined in an $A$ splitting and a pair weight of $\sqrt{B}$ to those joined in a $B$ splitting.  With this, the generalized Kauffman bracket of any link diagram $\tilde{L}$ is a specialization of the generalized transition polynomial.  The following theorem is a natural extension of Jaeger's~\cite{Jae90} relation between the original transition polynomial and Kauffman bracket.

\begin{thm}\label{L}
Let $L$ be a link
in $\Sigma\times [0,1]$ with link diagram $\tilde{L}$ and universe $\Gamma_L$. Then
$[\tilde{L}](A,B,d)=\frac{1}{d}Q(\Gamma_L; W_L, d)$ 

\end{thm}

\begin{proof}
This follows immediately from comparing Definitions~\ref{bracket} and \ref{transpolydef}.

\end{proof}

Chmutov and Pak~\cite{CP07} consider
\emph{signed ribbon graphs}, denoted $\hat{G}$,
where the sign on the edges acts as a device to keep track of the
over/under crossings of a not necessarily alternating link diagram. These are not the signed edges used to encode
topological information in unoriented ribbon graphs
defined previously. In this context, all ribbon graphs are oriented as topological surfaces.  However, each edge has a +/- indicator associated with it.
Chmutov and Pak~\cite{CP07} extend Definition~\ref{gen funct TT} to these signed ribbon graphs as follows.

\begin{defi} Let $\hat{G}$  be a signed ribbon graph, and let $F_1$ be the number of negative edges in $F$ and $F_2$ be the number of negative edges in $E(G)-F$.  Then

$$R(\hat{G};x,y,z)=\sum_{F \subseteq E(G)} (x-1)^{r(G)-r((F)+s((F)}y^{n((F)-s((F)}
z^{k((F)-bc((F)+n((F)},$$

where 
$s(F)=\frac{1}{2}(F_1 - F_2)$.

\end{defi}

Thus, if every edge of an signed oriented ribbon graph $\hat{G}$ is positive, then $R(\hat{G};x,y,z)=R(G; x, y, z, 1)$ where $G$ is the underlying unsigned ribbon graph, and $R(G; x, y, z, 1)$ is as in Definition~\ref{gen funct TT}.

If $\hat{G}$ is a signed oriented ribbon graph with medial
graph $G_m$, then we can define a signed weight system $W^-$ by reversing the roles of $\alpha$ and $\beta$ at vertices of $G_m$ corresponding to negative edges of $\hat{G}$.

With this, we have the following signed analog to Theorem~\ref{transpoly}, with a virtually identical proof that we leave to the reader.

\begin{thm}\label{signed trans}

Let $\hat{G}$ be a signed oriented ribbon graph with medial
graph $G_m$. Then

$$Q(G_m;W^-,t)=A^{r(G)}B^{n(G)}t^{k(G)}
R_{\hat{G}}(\frac{Bt}{A}+1,\frac{At}{B},\frac{1}{t}).$$

\end{thm}

\begin{prop}\label{signed medial universe}
If $\tilde{L}$ is a checkerboard colorable link diagram in an oriented surface $\Sigma$, with universe $\Gamma_L$, then $\Gamma_L$ is the medial graph ${G_L}_m$ for some $G_L$.  Furthermore, the link diagram induced weight system $W_L$ of $\Gamma_L$ is precisely the signed medial weight system $W^-$ of ${G_L}_m$ with $\alpha$ replaced by $A$ and $\beta$ replaced by $B$.
\end{prop}

\begin{proof}
The checkerboard coloring is exactly a face 2-coloring, say green and white, of the universe $\Gamma_L$  in $\Sigma$, with the green faces bounded by the half edges that are joined by an $A$ splitting.  Thus, $\Gamma_L$  is the medial graph of the green-face graph, which we denote $G_L$. We then note that the link diagram induced weight system $W_L$ of $\Gamma_L$ is precisely the medial weight system $W$ of ${G_L}_m$ with $\alpha$ replaced by $A$ and $\beta$ replaced by $B$, since the half edges paired by an $A$  splitting(respectively $B$ splitting) are precisely those with pair weight $\alpha$ (respectively $\beta$). 
\end{proof}

The medial graph construction of Proposition~\ref{signed medial universe}  provides a natural interpretation of the gluing procedure used by Chmutov and Pak~\cite{CP07} to produce $G_L$.

An alternative proof for the main theorem of Chmutov and Pak~\cite{CP07} then follows immediately.

\begin{thm}[Chmutov and Pak, \cite{CP07}]

If $\tilde{L}$ is checkerboard colored link diagram with signed green-face graph $\hat{G_L}$, then

$$[\tilde{L}](A,B,d)=A^{r(\hat{G}_L)}B^{n(\hat{G}_L)}d^{k(\hat{G}_L)-1}
R_{\hat{G}_L}(\frac{Bd}{A}+1,\frac{Ad}{B},\frac{1}{d}).$$

\end{thm}

\begin{proof} This is a consequence of Theorems~\ref{L} and~\ref{signed trans} and Proposition~\ref{signed medial universe}. 
\end{proof}

\begin{center}\textbf{Dedication}

\emph{This paper is dedicated to Tom Brylawski, who had a major influence on the first author, beginning in graduate school and continuing throughout her professional life.  He has given us a profound and lasting legacy of deep mathematics.}

\end{center}

\begin{center}\textbf{Acknowledgements}

We are grateful to Dan Archdeacon, Iain Moffat, and Lorenzo Traldi for a number of useful and interesting conversations.

\end{center}

\end{document}